\documentclass[a4paper,twoside,10pt]{article}

\usepackage[a4paper,width=150mm,top=25mm,bottom=25mm,bindingoffset=6mm]{geometry}

\usepackage[english]{babel}

\usepackage{graphicx}
\usepackage[font=small,labelfont=bf]{caption}
\usepackage{subcaption}
\usepackage{color}
\usepackage{mdframed}
\usepackage{framed} 
\usepackage{cancel}
\usepackage[normalem]{ulem}

\usepackage{lineno,hyperref}
\modulolinenumbers[5]

\usepackage{amsmath,amssymb,amsthm,mathtools,mathrsfs,epsfig,amsfonts,wasysym}
\usepackage{MnSymbol}

\usepackage[normalem]{ulem}
\usepackage{ragged2e}
\usepackage{csquotes}
\usepackage{tasks}
\usepackage{paralist}
\usepackage[inline]{enumitem}
\usepackage{tikz-cd}
\usepackage{stmaryrd}

\usepackage[titletoc,title]{appendix}

\newcommand{\e}{\mathrm{e}}

\newcommand{\overbar}[1]{\mkern 1.5mu\overline{\mkern-1.5mu#1\mkern-1.5mu}\mkern 1.5mu}

\usepackage{url,hyperref}
\hypersetup{colorlinks=true,linkcolor=blue,filecolor=magenta,urlcolor=cyan}
\usepackage{nicematrix}
\usepackage{tikz}
\usetikzlibrary{fit}

\usepackage[all]{xy}

\newtheorem{definition}{Defintion}
\newtheorem{theorem}{Theorem}
\newtheorem{proposition}{Proposition}

\usepackage{authblk}
\usepackage{todonotes}

\definecolor{wildstrawberry}{rgb}{1.0, 0.26, 0.64}

\definecolor{ao(english)}{rgb}{0.0, 0.5, 0.0}

\usepackage{orcidlink}

\begin{document}

\title{An inequality for completely monotone functions}

\author[1]{Vasiliki Bitsouni\,\orcidlink{0000-0002-0684-0583}\thanks{\texttt{vbitsouni@math.uoa.gr}}}

\author[1,2]{Nikolaos Gialelis\,\orcidlink{0000-0002-6465-7242}\thanks{\texttt{ngialelis@math.uoa.gr}}}

\author[3]{Dan-Stefan Marinescu\,\orcidlink{0000-0001-9066-2558}\thanks{\texttt{marinescuds@gmail.com}}}

\affil[1]{Department of Mathematics, National and Kapodistrian University of Athens, Panepistimioupolis, GR-15784 Athens, Greece}

\affil[2]{School of Medicine, National and Kapodistrian University of Athens,\newline GR-11527 Athens, Greece}

\affil[3]{National College \textquote{Iancu de Hunedoara}, Hunedoara, Romania}

\date{}

\maketitle

\begin{abstract}
\noindent
An inequality, which combines the concept of completely monotone functions with the theory of divided differences, is proposed. It is a straightforward generalization of a result, recently introduced by two of the present authors. 
\end{abstract}


\noindent
\textbf{MSC2020-Mathematics Subject Classification System:} 26A48, 26D07 



\section{Introduction}
\label{intro}

In \cite{bitsouni2022note}, the study of the generalized one-dimensional problem of Population Ecology led to the following multivariate analogue of a basic inequality,  
\begin{equation}
\label{gen}
\begin{gathered}
\prod\limits_{i=1}^n{{\left(1+x_i\right)}^{a_i}}\leq\e^{\frac{1}{n}\prod\limits_{i=1}^n{x_i}},\text{ where}\\
x_1,\dots,x_n\text{ are pairwise distinct non-negative real numbers and}\\ a_i\coloneqq\frac{\prod\limits_{\substack{j=1\\j\neq i}}^{n}{x_j}}{\prod\limits_{\substack{j=1\\j\neq i}}^{n}{\left(x_j-x_i\right)}}\,,\\
\text{with the equality holding only when one $x_i$ equals zero}.
\end{gathered}
\end{equation}
Here we show that \eqref{gen} is a special case of an inequality which involves a certain class of smooth non-negative functions, i.e. the completely monotone ones (see \hyperref[cmf]{Definition \ref*{cmf}}). In particular, we show the following result. 

\begin{proposition}
\label{crrs}
Let $f\colon\,\left(0,\infty\right)\to\left[0,\infty\right)$ be completely increasing (decreasing). 
\begin{enumerate}
\item Then 
\begin{equation}
\label{crrseqq}
f^{\left(n-1\right)}{\left(0\right)}\coloneqq\lim\limits_{x\to 0^+}{f^{\left(n-1\right)}{\left(x\right)}}\,\in\,\overbar{\mathbb{R}}\,\ni\,\lim\limits_{x\to \infty}{f^{\left(n-1\right)}{\left(x\right)}}\eqqcolon f^{\left(n-1\right)}{\left(\infty\right)},\text{ }\forall n\in\mathbb{N}.
\end{equation}
\item Let, also, $x_1,\dots,x_n\,\in\left(0,\infty\right)$ be pairwise distinct, with $$m\coloneqq \min\limits_{i\in\left\{1,\dots,n\right\}}{\left\{x_i\right\}}\text{ and }M\coloneqq\max\limits_{i\in\left\{1,\dots,n\right\}}{\left\{x_i\right\}}.$$ Then 
\begin{equation}
\label{crrseq}
\begin{gathered}
\frac{{\left(-1\right)}^{n-1}}{\left(n-1\right)!}f^{\left(n-1\right)}{\left(a\right)}\,\leq\,\sum\limits_{i=1}^n{\frac{f{\left(x_i\right)}}{\prod\limits_{\substack{j=1\\j\neq i}}^{n}{\left(x_j-x_i\right)}}}\,\leq\,\frac{{\left(-1\right)}^{n-1}}{\left(n-1\right)!}f^{\left(n-1\right)}{\left(b\right)},\text{ }\forall \left(a,b\right)\,\in\,\left[0,m\right]\times\left[M,\infty\right]\\
\left(\frac{{\left(-1\right)}^{n-1}}{\left(n-1\right)!}f^{\left(n-1\right)}{\left(b\right)}\,\leq\,\sum\limits_{i=1}^n{\frac{f{\left(x_i\right)}}{\prod\limits_{\substack{j=1\\j\neq i}}^{n}{\left(x_j-x_i\right)}}}\,\leq\,\frac{{\left(-1\right)}^{n-1}}{\left(n-1\right)!}f^{\left(n-1\right)}{\left(a\right)}\right).
\end{gathered}
\end{equation}
If, in addition, $f$ is strictly completely increasing (decreasing), then the inequalities in \eqref{crrseq} are strict.  
\end{enumerate}
\end{proposition}

The present short note is organized as follows: In \hyperref[Preliminaries]{Section \ref*{Preliminaries}} we introduce some basic notions and we state a preliminary result, which both are necessary for the statement and the proof of \hyperref[crrs]{Proposition \ref*{crrs}} in \hyperref[Proof of the main result]{Section \ref*{Proof of the main result}}. In \hyperref[Examples and corollaries]{Section \ref*{Examples and corollaries}} we give examples of the proposed inequality, one of which is \eqref{gen}.

\section{Preliminaries}
\label{Preliminaries}

First we give the definition of completely monotone functions. 
\begin{definition}
\label{cmf}
\begin{enumerate}
\item A function $f\in C^\infty{\left(\left(0,\infty\right);\left[0,\infty\right)\right)}$ is completely increasing (decreasing) iff $${\left(-1\right)}^n f^{\left(n\right)}{\left(x\right)}\leq 0\text{ }\left(\,{\left(-1\right)}^n f^{\left(n\right)}{\left(x\right)}\geq 0\,\right),\text{ }\forall \left(x,n\right)\,\in\,\left(0,\infty\right)\times\mathbb{N}.$$ Moreover, $f$ is strictly completely increasing (decreasing) iff the above inequality is strict.
\item A function $f\in C^\infty{\left(\left(0,\infty\right);\left[0,\infty\right)\right)}$ is (strictly) completely monotone iff it is either (strictly) completely increasing or (strictly) completely decreasing.
\end{enumerate}
\end{definition}
\noindent We note that the notion of complete monotonicity as suggested in \hyperref[cmf]{Definition \ref*{cmf}} does not appear in bibliography, where completely increasing functions are called Bernstein functions (see, e.g., \cite{schilling2012bernstein}), and completely decreasing functions are called completely (totally/absolutely) monotone (monotonic) functions (see, e.g., \cite{miller2001completely} and \cite{schilling2012bernstein}). To the authors' knowledge, the notion of strict complete monotonicity is also new. 

Additionally, we state a generalization of a well known result to higher derivatives, the mean value theorem for divided differences (see, e.g., \cite{popoviciu1933quelques}, \cite{sahoo1998mean}, or \cite{abel2004mean}). 
\begin{theorem}
\label{cauchypop}
Let $x_1,\dots,x_n\,\in\mathbb{R}$ be pairwise distinct, with $$m\coloneqq \min\limits_{i\in\left\{1,\dots,n\right\}}{\left\{x_i\right\}}\text{ and }M\coloneqq\max\limits_{i\in\left\{1,\dots,n\right\}}{\left\{x_i\right\}},$$ as well as $f\in C{\left(\left[m,M\right];\mathbb{R}\right)\cap C^{n-1}{\left(\left(m,M\right);\mathbb{R}\right)}}$. Then $\exists x_0\in\left(m,M\right)$, such that $$\left[x_1,\dots,x_n;f\right]\coloneqq\sum\limits_{i=1}^n{\frac{f{\left(x_i\right)}}{\prod\limits_{\substack{j=1\\j\neq i}}^{n}{\left(x_i-x_j\right)}}}=\frac{f^{\left(n-1\right)}{\left(x_0\right)}}{\left(n-1\right)!}.$$
\end{theorem}

\section{Proof of the main result}
\label{Proof of the main result}

We proceed by proving \hyperref[crrs]{Proposition \ref*{crrs}}. 
\begin{proof}
\begin{enumerate}
\item Let $n\in\mathbb{N}$ be abstract. We define $g\coloneqq{\left(-1\right)}^{n} f^{\left(n-1\right)}$ and we have $$g'{\left(x\right)}={\left(-1\right)}^{n} f^{\left(n\right)}{\left(x\right)}\leq 0\text{ }\left(\,g'{\left(x\right)}\geq 0\,\right),\text{ }\forall x\in\left(0,\infty\right),$$ hence $g$ is decreasing (increasing), which means that $$\lim\limits_{x\to 0^+}{g{\left(x\right)}}\,\in\,\overbar{\mathbb{R}}\,\ni\,\lim\limits_{x\to \infty}{g{\left(x\right)}}$$ and \eqref{crrseqq} then follows. 
\item From \hyperref[cauchypop]{Theorem \ref*{cauchypop}}, $\exists x_0\in \left(m,M\right)$ such that  $$-\sum\limits_{i=1}^n{\frac{f{\left(x_i\right)}}{\prod\limits_{\substack{j=1\\j\neq i}}^{n}{\left(x_j-x_i\right)}}}={\left(-1\right)}^{n}\left[x_1,\dots,x_n;f\right]=\frac{{\left(-1\right)}^{n}f^{\left(n-1\right)}{\left(x_0\right)}}{{\left(n-1\right)}!}.$$ In the light of point $1.$, we then get $$\begin{gathered}
\frac{{\left(-1\right)}^{n}}{\left(n-1\right)!}f^{\left(n-1\right)}{\left(b\right)}\,\leq\,-\sum\limits_{i=1}^n{\frac{f{\left(x_i\right)}}{\prod\limits_{\substack{j=1\\j\neq i}}^{n}{\left(x_j-x_i\right)}}}\,\leq\,\frac{{\left(-1\right)}^{n}}{\left(n-1\right)!}f^{\left(n-1\right)}{\left(a\right)},\text{ }\forall \left(a,b\right)\,\in\,\left[0,m\right]\times\left[M,\infty\right]\\
\left(\frac{{\left(-1\right)}^{n}}{\left(n-1\right)!}f^{\left(n-1\right)}{\left(a\right)}\,\leq\,-\sum\limits_{i=1}^n{\frac{f{\left(x_i\right)}}{\prod\limits_{\substack{j=1\\j\neq i}}^{n}{\left(x_j-x_i\right)}}}\,\leq\,\frac{{\left(-1\right)}^{n}}{\left(n-1\right)!}f^{\left(n-1\right)}{\left(b\right)}\right)
\end{gathered}$$ and \eqref{crrseq} then follows. 
\end{enumerate}
\end{proof}

\section{Examples and corollaries}
\label{Examples and corollaries}

The first example is that of the positive constant functions, which of course are both completely increasing and decreasing. Therefore, considering $f\in C^\infty{\left(\left(0,\infty\right);\left(0,\infty\right)\right)}$ where $$f{\left(x\right)}\coloneqq 1,\text{ }\forall x\in\left(0,\infty\right),$$ we have $$\lim\limits_{x\to 0^+}f^{\left(n\right)}{\left(x\right)}=0=\lim\limits_{x\to \infty}f^{\left(n\right)}{\left(x\right)},\text{ }\forall n\in\mathbb{N}$$ and we deduce from \hyperref[crrs]{Proposition \ref*{crrs}} that $$\sum\limits_{i=1}^n{\frac{1}{\prod\limits_{\substack{j=1\\j\neq i}}^{n}{\left(x_j-x_i\right)}}}=0\,,\text{ }\forall n\in\mathbb{N}\setminus\left\{1\right\}.$$

We now pass to more complicated examples. For this purpose, we have to refer to a catalogue of completely monotone functions. Extended lists of completely increasing and decreasing functions are contained in \cite{schilling2012bernstein} and \cite{miller2001completely}, respectively. For example, the function $f\in C^\infty{\left(\left(0,\infty\right);\left(0,\infty\right)\right)}$ where $$f{\left(x\right)}\coloneqq \frac{\ln{\left(1+x\right)}}{x},\text{ }\forall x\in\left(0,\infty\right),$$ is completely decreasing, and in particular, strictly completely decreasing. Besides, we have 
\begin{equation}
\label{bglm}
\lim\limits_{x\to 0^+}f^{\left(n-1\right)}{\left(x\right)}=\frac{{\left(-1\right)}^n\left(n-1\right)!}{n}\text{ and }\lim\limits_{x\to \infty}f^{\left(n-1\right)}{\left(x\right)}=0,\text{ }\forall n\in\mathbb{N}.
\end{equation}
Therefore, from \hyperref[crrs]{Proposition \ref*{crrs}} we deduce $$0\,<\,\sum\limits_{i=1}^n{\frac{f{\left(x_i\right)}}{\prod\limits_{\substack{j=1\\j\neq i}}^{n}{\left(x_j-x_i\right)}}}\,<\,\frac{1}{n}\,,\text{ }\forall n\in\mathbb{N}$$ and \eqref{gen} then follows. 

To sum up, the main result proposes an elegant, systematic and unified approach for the proof of \eqref{gen} and relevant inequalities. The price we pay is the evaluation of quantities such as \eqref{bglm}. 

Below follow more examples. We note that the choice of these examples is based on 
\begin{itemize}
\item the finite behavior at $0$ and $\infty$ of the corresponding completely monotone functions and their derivatives, in order to avoid a trivial right-hand side inequality of the following type $\left[0,\infty\right)\ni A\leq B\leq\infty$, as well as
\item the simplicity of the evaluation of the right/left-hand side of \eqref{crrseq}, in order to keep the presentation as compact as possible.
\end{itemize}
\begin{enumerate}
\item \fbox{Strictly completely increasing functions}:
\begin{enumerate}[label=\roman*.]
\item Let $\alpha\in\left(0,\infty\right)$ and $f\in C^\infty{\left(\left(0,\infty\right);\left(0,\infty\right)\right)}$ where $$f{\left(x\right)}\coloneqq \frac{x}{\alpha+x},\text{ }\forall x\in\left(0,\infty\right).$$ Then $$\lim\limits_{x\to 0^+}f^{\left(n\right)}{\left(x\right)}=\frac{{\left(-1\right)}^{n+1}n!}{\alpha^{n}}\text{ and }\lim\limits_{x\to \infty}f^{\left(n\right)}{\left(x\right)}=0,\text{ }\forall n\in\mathbb{N},$$ hence 
$$0\,<\,-\sum\limits_{i=1}^n{\frac{f{\left(x_i\right)}}{\prod\limits_{\substack{j=1\\j\neq i}}^{n}{\left(x_j-x_i\right)}}}\,<\,\frac{1}{\alpha^{n-1}},\text{ }\forall n\in\mathbb{N}\setminus\left\{1\right\}.$$ 
\item Let $\beta\in\left(0,\infty\right)$, $\alpha\in\left(0,\beta\right)$ and $f\in C^\infty{\left(\left(0,\infty\right);\left(0,\infty\right)\right)}$ where $$f{\left(x\right)}\coloneqq \ln{\left(\frac{\beta\left(x+\alpha\right)}{\alpha\left(x+\beta\right)}\right)},\text{ }\forall x\in\left(0,\infty\right).$$ Then $$\lim\limits_{x\to 0^+}f^{\left(n\right)}{\left(x\right)}=\frac{{\left(-1\right)}^{n+1}n!}{n}\left(\frac{1}{\alpha^n}-\frac{1}{\beta^n}\right)\text{ and }\lim\limits_{x\to \infty}f^{\left(n-1\right)}{\left(x\right)}=0,\text{ }\forall n\in\mathbb{N},$$ hence $$0\,<\,-\sum\limits_{i=1}^n{\frac{f{\left(x_i\right)}}{\prod\limits_{\substack{j=1\\j\neq i}}^{n}{\left(x_j-x_i\right)}}}\,<\,\frac{1}{n-1}\left(\frac{1}{\alpha^{n-1}}-\frac{1}{\beta^{n-1}}\right),\text{ }\forall n\in\mathbb{N}\setminus\left\{1\right\}.$$ 
\end{enumerate}
\item \fbox{Strictly completely decreasing functions}:
\begin{enumerate}[label=\roman*.]
\item Let $\alpha\in\left(0,\infty\right)$ and $f\in C^\infty{\left(\left(0,\infty\right);\left(0,\infty\right)\right)}$ where $$f{\left(x\right)}\coloneqq \e^{-\alpha x},\text{ }\forall x\in\left(0,\infty\right).$$ Then $$\lim\limits_{x\to 0^+}f^{\left(n-1\right)}{\left(x\right)}={\left(-1\right)}^{n-1}\alpha^{n-1}\text{ and }\lim\limits_{x\to \infty}f^{\left(n-1\right)}{\left(x\right)}=0,\text{ }\forall n\in\mathbb{N},$$ hence $$0\,<\,\sum\limits_{i=1}^n{\frac{f{\left(x_i\right)}}{\prod\limits_{\substack{j=1\\j\neq i}}^{n}{\left(x_j-x_i\right)}}}\,<\,\frac{\alpha^{n-1}}{\left(n-1\right)!}\,,\text{ }\forall n\in\mathbb{N}.$$ 
\item Let $\alpha,\beta,\gamma\,\in\left(0,\infty\right)$ and $f\in C^\infty{\left(\left(0,\infty\right);\left(0,\infty\right)\right)}$ where $$f{\left(x\right)}\coloneqq {\left(\alpha+\beta x\right)}^{-\gamma},\text{ }\forall x\in\left(0,\infty\right).$$ Then $$\lim\limits_{x\to 0^+}f^{\left(n-1\right)}{\left(x\right)}=\alpha^{\gamma+n-1}\beta^{n-1}\binom{-\gamma}{n-1}\left(n-1\right)!\text{ and }\lim\limits_{x\to \infty}f^{\left(n-1\right)}{\left(x\right)}=0,\text{ }\forall n\in\mathbb{N},$$ hence $$0\,<\,\sum\limits_{i=1}^n{\frac{f{\left(x_i\right)}}{\prod\limits_{\substack{j=1\\j\neq i}}^{n}{\left(x_j-x_i\right)}}}\,<\,\alpha^{\gamma+n-1}\beta^{n-1}\binom{\gamma+n-2}{n-1}\,,\text{ }\forall n\in\mathbb{N}.$$ 
\end{enumerate}
\end{enumerate}

 \bibliographystyle{plain}
\bibliography{mybibfile}\label{bibliography}
\addcontentsline{toc}{chapter}{Bibliography}

\end{document}